[1994/12/01]

\documentclass[draft]{article}
\usepackage{amsmath, amsthm, amsfonts}
\title{\textbf{Algorithmic definition of means acting on positive numbers and operators}}
\author{Mikl\'os P\'alfia\thanks{I would like to thank D\'enes Petz and N\'andor Sieben for their help and support.}}
\date{March 6, 2005}

\theoremstyle{plain} 
\newtheorem{thm}{Theorem}[section]
\newtheorem{cor}[thm]{Corollary}

\theoremstyle{definition}
\newtheorem{defn}{Definition}[section]

\theoremstyle{remark}
\newtheorem{rem}{Remark}[section]

\numberwithin{equation}{section}

\begin{document}

\maketitle

\begin{abstract}
Means are used in several applications from electronic engeneering to information theory, however there is no general theorem on how to extend a given $M(x,y)$ mean function to multiple variable forms. In this article we would like to present a theorem, which gives one possible solution for this problem, for every $M(x,y)$ mean function, acting on positive numbers and operators.
\end{abstract}

\section{Introduction}
The use of \textit{mean} functions falls back to the early ages of mathematics, but they are used widely nowdays aswell. For example they have great importance in statistics, but also in volume visualization, in imaging surgery.

The study of electrical network connections implied the introduction of parallel sum of two positive semidefinite matrices in \cite{anderson_duffin}. Formerly, Anderson defined a matrix operation, called shorted operation to a subspace, for each positive semidefinite matrix. Anderson and Trapp in \cite{anderson_trapp}, have extended the theory of parallel addition and shorted operation to bounded linear positive operators on a Hilbert space and demonstrated
its importance in operator theory. They have studied fundamental properties of these operations.

The axiomatic theory of \textit{means}, for pair of positive operators, have been developed by Kubo and Ando in \cite{kubo}. This theory has found a number of applications in operator theory.

Several steps have been taken in the theory of \textit{means}, but a general idea has not been laid down yet, on how to extend, or define an $n$ variable version of a given general $M(x,y)$ \textit{mean} function acting on numbers and operators. We would like to present a theory which gives a possible solution and a frame theory for further studies.
\begin{defn}
A two variable function \textit{M} :
$\mathbb{R}^+\times\mathbb{R}^+\mapsto\mathbb{R}^+$ according to \cite{petz}, is called
a \textit{mean} function if
\begin{enumerate}
\renewcommand{\labelenumi}{(\roman{enumi})}
\item $M(x,x)=x$ for every $x\in\mathbb{R}^+$.
\item $M(x,y)=M(y,x)$ for every $x,y\in\mathbb{R}^+$.
\item If $x<y$, then $x<M(x,y)<y$.
\item If $x<x'$ and $y<y'$, then $M(x,y)<M(x',y')$.
\item $M(x,y)$ is continuous.
\end{enumerate}
\end{defn}

The \textit{geometric mean} $\sqrt{xy}$, the \textit{arithmetic
mean} $(x+y)/2$ and the \textit{harmonic mean} $2/(x^{-1}+y^{-1})$
are the most known examples but there are many other \textit{means} aswell.

The above definition of means can be easily extended to positive
ordered operators and matrices, by replacing the numbers with them.
\section{Extending a mean to multiple variables}
The definition of $M(x,y)$ given in the first section, can be extended to multiple variable functions. The following definition is one possibility of extension.
\begin{defn}
An $n$ variable function $M_n : (\mathbb{R}^+)^{\textit{n}}
\mapsto\mathbb{R}^+$ may be called a ($n$ variable) \textit{mean}
function if
\begin{enumerate}
\renewcommand{\labelenumi}{(\roman{enumi}')}
\item $M_n(x,\dots x)=x$ for every $x\in\mathbb{R}^+$.
\item $M_n$ is independent from the ordering of $x_1,x_2,\dots x_n$.
\item $min(x_1,\dots x_n)\leq M_n(x_1,\dots x_n)\leq max(x_1,\dots x_n)$.
\item If $\exists (i)\quad x_i\leq x'_i$, then $M_n(x_1,\dots x_n)\leq M_n(x_1,\dots x'_i,\dots x_n)$.
\item If $\forall (i)\quad x_i<x'_i$, then $M_n(x_1,\dots x_n)<M_n(x'_1,\dots x'_n)$.
\item $M_n(x_1,\dots x_n)$ is continuous.
\end{enumerate}
\end{defn}

Our goal is to algorithmically define an $n$ variable
\textit{mean} by using its less than $n$ variable forms.
Firstly we will use the $n-1$ variable form of $M$, to define the $n$ variable one, as an iteration's limit.
\begin{defn}\label{defn0}
Let $X=(x_1,\dots x_n)\in(\mathbb{R}^+)^{\text{n}}$ and $M_{n-1}$
be an $n-1$ variable mean function. Let us consider the $n-1$
class variations of $x_n$. There are $\binom{n}{n-1}=n$
different variations. Let us define an iteration as
\begin{gather}
    x_i^0=x_i \quad \forall i\in[1,\dots n]\\
    \label{xk}x_i^{k+1}=M_{n-1}(V_{i}^{n-1}(X^k))
\end{gather}
where $V_{i}^{n-1}(X^k)=V_{i}^{n-1}(x_1^k,\dots x_n^k)$ is the
$i$th $n-1$ different variation of $x_n$ ($n-1$ are chosen from
$n$, which could be done in $n$ different ways).
\end{defn}
\begin{thm}\label{thm1}
The $x_i^k$ sequences - defined in definition \ref{defn0} - are convergent and their limits are the same, which could
be defined as the $M_n$ mean of the $x_n$ numbers.
\end{thm}
\begin{proof}
The iteration given in definition \ref{defn0} has a contractive-like property by
(iii'), which means that the sequence $\min{X^k}$ is monotonic
increasing, $\max{X^k}$ is monotonic decreasing. Hence the
limits $\lim_{k \to \infty}\min{X^k}$ and $\lim_{k \to
\infty}\max{X^k}$ exist. From condition (iv'), one can see that
the series' minimal and maximal elements are given by the
following:
\begin{gather}
    \label{minxk}\min{X^{k+1}}=M_{n-1}(V_{min}^{n-1}(X^k))\\
    \label{maxxk}\max{X^{k+1}}=M_{n-1}(V_{max}^{n-1}(X^k))
\end{gather}
where $V_{min}^{n-1}(X^k)$ are the smallest $n-1$ numbers,
$V_{max}^{n-1}(X^k)$ are the largest $n-1$ numbers from the series
$X^k$. By \eqref{minxk} and \eqref{maxxk} $\min{X^{k+1}}$ and
$\max{X^{k+1}}$ explicit dependence on $\min{X^{k}}$ and
$\max{X^{k}}$ is given
\begin{gather}
    \label{limminxk}\lim_{k \to \infty}\min{X^{k+1}}=\lim_{k \to \infty}M_{n-1}(V_{min}^{n-1}(X^k))\\
    \label{limmaxxk}\lim_{k \to \infty}\max{X^{k+1}}=\lim_{k \to \infty}M_{n-1}(V_{max}^{n-1}(X^k))
\end{gather}
and by (vi') one can write
\begin{gather}
    \lim_{k \to \infty}\min{X^{k+1}}=M_{n-1}(\lim_{k \to \infty}V_{min}^{n-1}(X^k))\\
    \lim_{k \to \infty}\max{X^{k+1}}=M_{n-1}(\lim_{k \to \infty}V_{max}^{n-1}(X^k))
\end{gather}
which yields
\begin{gather}
    \lim_{k \to \infty}\min{X^{k+1}}=I=M_{n-1}(I,\dots)\\
    \lim_{k \to \infty}\max{X^{k+1}}=J=M_{n-1}(J,\dots)
\end{gather}
and this is only true when $I=J$.
\end{proof}
The above theorem and its proof yields the following remarks.
\begin{rem}
The iteratively defined \textit{mean} function $M_n$, is invariant
to the initial ordering of the $V_{i}^{n-1}$ variations.
\end{rem}
\begin{rem}
The iteration \eqref{xk} leaves the \textit{mean} of the starting
n numbers invariant through the sequence.
\end{rem}
It is easy to verify the next two theorems, which have stressed importance in inequalities of \textit{means} and \textit{operator means} in our given context.
\begin{cor}
Two different, iteratively defined \textit{mean} function $M_{n,1}$ and $M_{n,2}$, is in the same relation as their two variable forms ($M_{2,1}<M_{2,2}$ implies $M_{n,1}<M_{n,2}$).
\end{cor}
\begin{cor}
Theorem \ref{thm1} and the proof also works for ordered positive
operators, acting on a $H$ Hilbert space, and for $k\times k$ matrices.
\end{cor}
We will show with some examples, that theorem \ref{thm1} gives a sufficient definition of the $n$ variable \textit{mean}.
\begin{cor}
Theorem \ref{thm1} applied on the $n-1$ variable
\textit{arithmetic}, \textit{geometric} and \textit{harmonic
mean}, gives the corresponding $n$ variable \textit{mean}.
\end{cor}
\begin{proof}
According to the given serie's convergence in theorem \ref{thm1}
it is enough to prove that the minimum's or maximum's limit is the
$n$ variable \textit{mean}. We will prove it only for the
\textit{arithmetic mean}. Let us consider the numbers
$x_1^0\leq\dots \leq x_n^0 \in\mathbb{R}^+$ and the $A_{n-1}$
\textit{arithmetic mean}. By theorem 2.1, the sequence
$x_1^k$ can explicitly be written and proven by induction with \eqref{xk}, \eqref{minxk} and \eqref{maxxk}:
\begin{equation}
\begin{split}
x_1^k=
    \begin{cases}
    \frac{\sum_{i=1}^{n}\frac{(n-1)^k-1}{n}x_i^0+x_1^0}{(n-1)^k}&
    \text{if $k$ is even,}\\
    \frac{\sum_{i=1}^{n}\frac{(n-1)^k+1}{n}x_i^0-x_n^0}{(n-1)^k}&
    \text{if $k$ is odd,}
    \end{cases}\\
x_n^k=
    \begin{cases}
    \frac{\sum_{i=1}^{n}\frac{(n-1)^k-1}{n}x_i^0+x_n^0}{(n-1)^k}&
    \text{if $k$ is even,}\\
    \frac{\sum_{i=1}^{n}\frac{(n-1)^k+1}{n}x_i^0-x_1^0}{(n-1)^k}&
    \text{if $k$ is odd.}
    \end{cases}
\end{split}
\end{equation}

For the \textit{geometric mean} the proof can be extended using
the logarithmic function and its inverse for the limit:
\begin{equation}
\log(\mspace{-8.0mu}\sqrt[\leftroot{-2}\uproot{2}n-1]{x_1\cdot x_2\dots x_{n-1}})=\frac{\sum_{i=1}^{n-1}\log x_i}{n-1}\text{.}
\end{equation}
The proof for the \textit{harmonic mean} can be given by inverses:
\begin{equation}
\left(\frac{n-1}{x_1^{-1}+x_2^{-1}\dots x_{n-1}^{-1}}\right)^{-1}=\frac{x_1^{-1}+x_2^{-1}\dots x_{n-1}^{-1}}{n-1}\text{.}
\end{equation}
\end{proof}

Our main idea of extending \textit{means} to multiple variables is
based on theorem ~\ref{thm1}, which is theoretically enough but
in practice is very insufficient. For example if we would like to
compute $n$ numbers or matrices $M_n$ \textit{mean}, we should use
the two variable main definition $M_2$ and extend the other ones
from one to another. In the next section we will prove that $M_n$
can be directly extended from the corresponding $M_2$.
\section{Extending $M_n$ directly from $M_2$}
Let us define the following iteration:
\begin{defn}\label{defn1}
Let $X=(x_1^0\leq\dots \leq x_n^0)\in(\mathbb{R}^+)^{\text{n}}$
and $M=M_2$ be a two variable \textit{mean} function,
\begin{equation}
    x_i^{k+1}=
    \begin{cases}
    M(x_1^k,x_2^k)&\text{if $i=1$,}\\
    M(x_{n-1}^k,x_n^k)&\text{if $i=n$,}\\
    M(x_{i-1}^k,x_{i+1}^k)&\text{else.}
    \end{cases}
\end{equation}
\end{defn}
\begin{thm}\label{thm2}
The iteration given in definition \ref{defn1} for all $n$ is
convergent and
\begin{equation*}
\forall (i)\quad \lim_{k \to \infty}x_i^k=M_n(x_1^0,\dots x_n^0)
\end{equation*}
where $M_n$ is defined by theorem \ref{thm1}.
\end{thm}
\begin{proof}
Firstly we begin with proving the convergence.
It is clear that the $\min x_i^{k}$ is always the first element ($i=1$)
and the $\max x_i^{k}$ is always the last element ($i=n$) of the series
in definition \ref{defn1}. Hence (iii) and definition \ref{defn1},
$\min x_i^{k}$ is increasing and $\max x_i^{k}$ is decreasing. This
yields:
\begin{gather}
		\lim_{k \to \infty}\min x_i^{k+1}=\lim_{k \to \infty}M(x_1^k,x_2^k)\\
		\lim_{k \to \infty}\max x_i^{k+1}=\lim_{k \to \infty}M(x_{n-1}^k,x_n^k)
\end{gather}
and by (v):
\begin{gather}
		\lim_{k \to \infty}\min x_i^{k+1}=M(\lim_{k \to \infty}x_1^k,\lim_{k \to \infty}x_2^k)\\
		\lim_{k \to \infty}\max x_i^{k+1}=M(\lim_{k \to \infty}x_{n-1}^k,\lim_{k \to \infty}x_n^k)
\end{gather}
which give
\begin{gather}
		\lim_{k \to \infty}\min x_i^{k+1}=I=M(I,\lim_{k \to \infty}x_2^k),\\
		\lim_{k \to \infty}\max x_i^{k+1}=J=M(\lim_{k \to \infty}x_{n-1}^k,J).
\end{gather}
Considering the characteristics of the definition \ref{defn1}, this can only be true
when $I=J$.

Secondly we will prove the limit.
For $n=3$ the theorem is clear, because the two iterations, defined
in theorems \ref{thm1} and in \ref{thm2}, are the same. Our next step is
to prove for $n+1$, if it is true for $n$.

Let us consider the definition of $M_{n+1}$ in theorem \ref{thm1}. Comparing
the $\min x_i^{k}$ (which is the first element $i=1$) in theorem \ref{thm2},
and $\min X^k$ (which equals $M_{n}(V_{min}^{n}(X^{k-1}))$ by \eqref{minxk}),
we can see that $\min x_i^{k}\leq \min X^k$, because of the inductional
condition and the definition of $M_{n}$ as a limit in theorem \ref{thm1}. The same can
be applied for $\max x_i^{k}$ and $\max X^k$ which yields $\max x_i^{k}\geq \max X^k$.
Hence $\min x_i^{k}$ and $\max x_i^{k}$ are minoring and majoring, for every $k$,
$\min X^k$ and $\max X^k$, but $\lim_{k \to \infty}\min x_i^{k}=\lim_{k \to \infty}\max x_i^{k}$,
so $\lim_{k \to \infty}\min x_i^{k}=M_{n+1}(x_1^0,\dots x_{n+1}^0)$.
\end{proof}

Furthermore there is special property in the iteration in definition \ref{defn1}.
\begin{defn}\label{defn2}
Let $x_i^0\in\mathbb{R}^+\quad i\in[1,\dots n]$ and $G$ be a graph, with n verteces and edges given as that, there is one cycle in $G$, which contains all verteces and edges (so it is at the same time a \textit{Hamiltonian}- and an \textit{Euler-cycle}). This implies that in $G$, every vertex has two edges and all of them are bound together.
Let us consider an optional one to one correspondence between $x_i^0$ numbers and $G$-s verteces.
Taking every edge in $G$ as an $M(x_j^0,x_l^0)$ (where $M$ is a \textit{mean} function and $x_j^0$, $x_l^0$ are assigned to the two ending points of the edge as previously given), we can define an iteration with an optional $n$ mappings,
\begin{equation}
\begin{split}
x_i^{k+1}=M(x_j^k, x_l^k)\\
i,j,l\in [1,\dots n]\quad j\neq l\text{.}
\end{split}
\end{equation}
\end{defn}
\begin{thm}\label{thm3}
Every different iteration given in definition \ref{defn2}, converge to the limit $M_n$ \textit{mean} function,
and the iteration - independently from the mapping $x_i^{k+1}=M(x_j^k, x_l^k)$ - converge on a higher or equally rate as the iteration given in definition \ref{defn1}.
\end{thm}
\begin{proof}
For $n=3$, it is easy to see that the theorem is true, because the iterations given in definitions \ref{defn1} and \ref{defn2},
are the same.

Assume that the theorem is true for $n$ variable. Let us expand from an $n$ variable iteration defined in
\ref{defn2} with an optional mapping, to $n+1$ variable. This can be done as replacing one edge with
two edges and one vertex (mapped to a new number). Let us do this expansion in the following way. Take the first
smallest $n$ numbers from $n+1$ and set up on them the iteration given in definition \ref{defn1}. From the inductional condition this iteration will have the slowest convergence rate, which means that its minimal and maximal
elements will minor and major every other iteration in definition \ref{defn2}. Let us replace the edge which gives the maximal element of the iteration given in definition \ref{defn1}, with the two new edges and the remaining number (which is the greatest number out of the $n+1$) as a vertex. This two edges with the corresponding two $M_2(x,y)$ will give the new iterations - given for $n+1$ numbers - greatest two elements. This replacement cannot be done better in any other optional mapped iteration given in definition \ref{defn2} aswell. But considering the inductional condition this yields that any $n+1$ variable iteration given in definition \ref{defn2}, cannot minor and major, with its maximal and minimal elements, the iteration given in definition \ref{defn1}, hence theorem \ref{thm3} is proven.
\end{proof}
\begin{cor}
The above theorems also work for ordered operators acting on a $H$ Hilbert space, and $k\times k$ matrices aswell.
\end{cor}
We will have to consider further examinations to define the above iterational definitions for inorderable matrices. In the next section we will study this problem.

\section{Extending theorems \ref{thm2} and \ref{thm3} to unordered matrices and operators}
The problem is with positive matrices and operators which satisfy $\left\|A\right\|=\left\|B\right\|$ and $A\neq B$.
For the above matrices and operators, the function $M(A,B)$'s and its arguments' relation is not explained and highly depend on the main characteristics of $M(A,B)$, so the given iterations in the above theorems must be specified.
\begin{thm}\label{thm4}
For any $X_1,\dots X_n$ positive operators or matrices the iteration given in definition \ref{defn1} is convergent, as defined in theorem \ref{thm2}.
\end{thm}
\begin{proof}
If $\left\|X_1\right\|=\left\|X_2\right\|=\dots =\left\|X_n\right\|$ does not hold, than the iteration in definition \ref{defn1} converges for all $X_i$, because after $n$ steps, the iteration will surely alter all of the $\left\|X_i^k\right\|$-s (from one to another), so the iteration will converge.

If $\left\|X_1\right\|=\left\|X_2\right\|=\dots =\left\|X_n\right\|$ does hold, we will have to define (according to \cite{petz}) the following construction. Let $a_t$ and $a'_t$ be monotone sequences as,
\begin{equation}\label{eq3}
\begin{split}
&a_t, a'_t\in\mathbb{R}^+\\
&\forall (t)\quad a_t\geq 1 \text{ and } a'_t\leq 1\\
&lim_{t \to \infty}a_t=1 \text{ and }lim_{t \to \infty}a'_t=1\text{.}
\end{split}
\end{equation}
Let $X'_1=a_tX_1$, $X'_i=X_i\quad i\in[2,\dots n]$ and $X''_1=a'_tX_1$, $X''_i=X_i\quad i\in[2,\dots n]$. Let us set the iteration given in definition \ref{defn1} up on $X_i$, $X'_i$ and $X''_i$. According to the first part of the proof, the $(X'_i)^k$ and $(X''_i)^k$ series are convergent for any $t$, as given in theorem \ref{thm2}. Considering the definition of sequences $a_t$ and $a'_t$ in \eqref{eq3}, it is easy to verify by condition (iii), that for any $t$, the $(X''_i)^k$ series are minoring and $(X'_i)^k$ series are majoring the series $X_i^k$ for any $i$. Taking the limit $k\to \infty$, we get $M_n(X''_1(t),\dots X''_n(t))$ and $M_n(X'_1(t),\dots X'_n(t))$. Hence $M_n(X''_1(t),\dots X''_n(t))$ and $M_n(X'_1(t),\dots X'_n(t))$ are Cauchy sequences in index $t$ and condition (vi'), they are convergent and
\begin{equation}\label{eq4}
lim_{t \to \infty}M_n(X''_1(t),\dots X''_n(t))=lim_{t \to \infty}M_n(X'_1(t),\dots X'_n(t))\text{.}
\end{equation}
But $(X''_i)^k(t)$ and $(X'_i)^k(t)$ are minoring and majoring every $X_i^k$ for any $i$ and $t$, so the limit $lim_{k \to \infty}X_i^k$ exist and by \eqref{eq4},
\begin{equation}
lim_{t \to \infty}M_n(X''_1(t),\dots X''_n(t))=lim_{k \to \infty}X_i^k=lim_{t \to \infty}M_n(X'_1(t),\dots X'_n(t))
\end{equation}
and theorem \ref{thm4} is proven.
\end{proof}
\begin{cor}
Using the above proof, theorems \ref{thm1} and \ref{thm3} work for the unordered $X_i$-s.
\end{cor}
\section{Consequences}
By the theorems given in our examinations generalize the extension of the two variable \textit{mean} functions and gives a frame theory, which may be used in the future studies related to the extension of \textit{means} to multiple variables. An important outcome is, that these theorems are applying for operators and matrices and guarantee the existence of one possible extension.

It is known, that in several situations, there are more than one possible generalization of a \textit{mean}. One example is the \textit{logarithmic mean},
\begin{equation}
L(x,y)=\frac{x-y}{\log x-\log y}
\end{equation}
which has several extended forms according to \cite{carlson}, \cite{mustonen}, \cite{neuman}, but our theorems may leave only one form valid. However with some \textit{means}, it appears to be quite difficult to give the iterations limit in a closed form.

\end{document}